# The Surprising Non-Triviality of Sharing a Sandwich with Two Other People

John O'Meara
Montclair State University
Correspondence: omearaj1@montclair.edu

## Introduction

When I'm not wearing the hat of a Ph.D. student, adjunct instructor, and/or graduate researcher, I work through the summers running a sandwich stand at the county fair with my brother. I find this gig appealing because it presents an entirely new set of challenges and opportunities to problem solve that feed my desire to remain a well-rounded individual. We will often conduct research and development on new sandwich recipes in our test kitchen before debuting them to the public. On one such day, we were experimenting with variants of the BLT when we decided to share it amongst the both of us and our mom. My brother passively commented that it seemed difficult to cut the sandwich into perfect thirds when making diagonal cuts. I went over and looked at the sandwich, and felt compelled to agree with him. It *does* seem difficult to cut this rectangular sandwich into perfect thirds with cuts other than the trivial straight cuts that create rectangular strips! Despite having put my mathematician hat off to the side, I ran to put it back on and see what the math had to say about this dilemma.

## The Fair Sandwich Problem

*Given a rectangular sandwich of dimensions x and y, is it possible to split the sandwich into three pieces of equal area when restricting ourselves to diagonal cuts?*

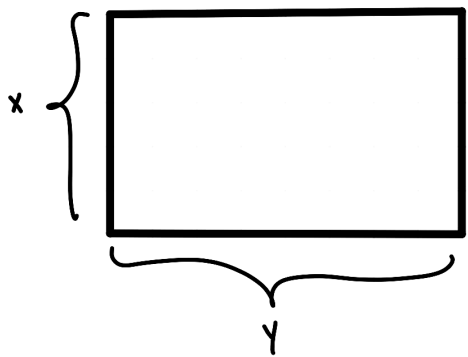

**Figure 1.** A rectangular sandwich of your choosing.

**Parameters of Problem-Solving**

**A simple diagonal cut.** We will approach this problem in a couple of different ways. The first assumption we will make is that our diagonal cuts will originate from the corner of the sandwich and extend in a straight line segment to the opposite side at some point (see Figure 2). We will assume that, without loss of generality, the distance of the two cuts need not be equal.

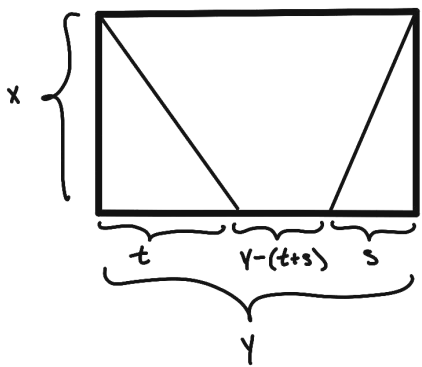

**Figure 2.** Imposing two diagonal cuts on our rectangular sandwich to yield three regions of sandwich to be shared. Note that, by construction, $t$ and $s$ must lie in the interval $(0, 0.5y]$.

Note that this method of cutting yields two right triangles and one trapezoid. In order for this method to yield three evenly cut pieces of sandwich, all three of these areas must be equal. Let us call these areas, from left to right in Figure 2, $A_t$, $A_{y-(t+s)}$, and $A_s$. Then, we have that:

$A_t = \frac{1}{2}xt$ , $A_{y-(t+s)} = \frac{1}{2}(2y - t - s)(x)$ , and $A_s = \frac{1}{2}xs$ .

Given that all three areas are equal, and thus $A_t = A_s$ , this implies $t = s$ . By division (assuming $x \neq 0$) and substitution:

$t = 2y - 2t \Rightarrow 3t = 2y \Rightarrow t = \frac{2}{3}y$ .

However, because $t \in (0, \frac{1}{2}y]$ by construction, this is impossible. Therefore, this cutting method cannot yield perfect thirds of the original rectangular sandwich. ◇

**A perturbation of the diagonal cut.** Let us keep going forward with the method of diagonal cuts. We will now assume that, although the cuts remain line segments that extend from one side to the opposite parallel side of the same length, that the cuts need not originate nor terminate in a corner of the sandwich (see Figure 3).

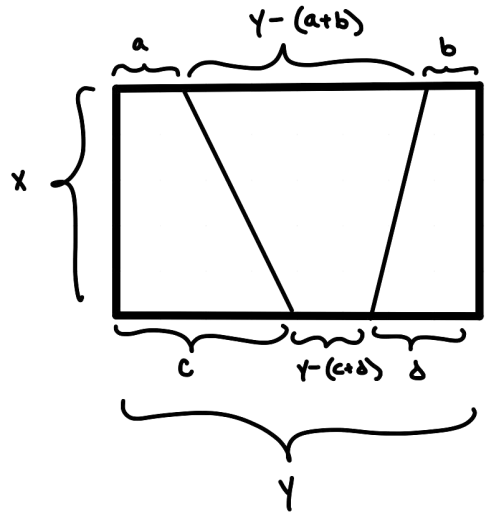

**Figure 3.** A perturbation of the diagonal cut method.

We see that this perturbation now yields three trapezoidal regions of sandwich that all need to have the same area. We will call these regions, from left to right, $A_c$, $A_{y-(c+d)}$, $A_d$. Thus, we have:

$$A_c = \frac{1}{2}(a+c)(x) \,, \; A_{y-(c+d)} = \frac{1}{2}(2y - (a+b+c+d)(x) \,, \text{ and } A_d = \frac{1}{2}(b+d)(x) \,.$$

Assuming these areas are equal, we obtain the fact that $a + c = b + d$. By substitution and division, we have:

$$b + d = 2y - 2(b+d) \Rightarrow 3(b+d) = 2y \Rightarrow b + d = a + c = \frac{2}{3}y \,.$$

Let's investigate this set of solutions. Certainly, when $a = b = c = d = \frac{1}{3}y$, we have a solution (the trivial rectangular strip cuts). However, are all solutions such that $b + d = a + c = \frac{2}{3}y$ feasible?

If $b + d = a + c = \frac{2}{3}y$ , this implies that $b = \frac{2}{3}y - d$ and $a = \frac{2}{3}y - c$ . Note that, for any nonzero choice of $a, b, c, d \neq \frac{1}{3}y$ , we obtain solutions for the dependent parameters ($b$ depending on $d$, $a$ depending on $c$) that are outside of the interval $(0, \frac{1}{3}y]$ . This means that if either $c$ or $d$ are greater than $\frac{1}{3}y$, then $a$ or $b$ must be less than $\frac{1}{3}y$ , with equality occurring at $a = b = c = d = \frac{1}{3}y$, which is a case that provides perfect sandwich thirds. Because a line that goes through the intersection of the diagonals always splits area into equal parts (see Figure 4), and this could occur by construction of the perturbed diagonal cut as the dependence of $b$ and $a$ on $d$ and $c$ respectively, this method ensures that perfect thirds can be achieved.

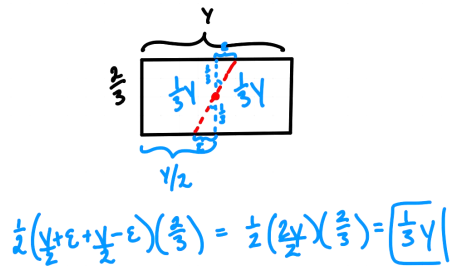

**Figure 4.** Visual proof that any line that goes through the intersection of the rectangle's diagonals splits the area of the rectangle exactly in half.

However, what if the cut does not go through the intersection of the diagonals? If such a cut is chosen, the areas are no longer split evenly. As seen in Figure 5, it is impossible for the similar triangles formed between the red and blue dotted segments to be congruent, and thus the areas of the two regions cannot be equal.

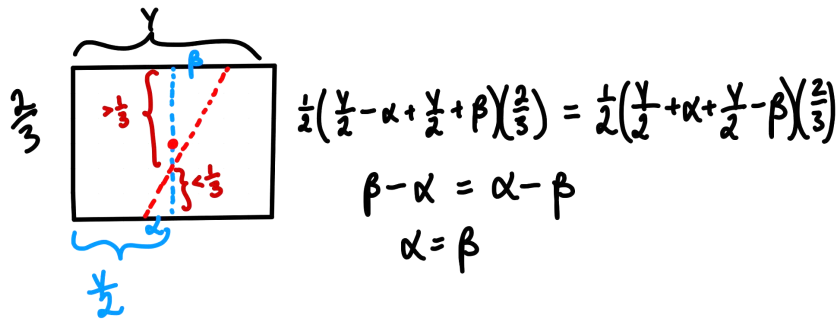

**Figure 5.** A visual proof to show that a segment that does not go through the intersection of a rectangle's diagonals cannot split the area in two.

There is the possibility that an area of $\frac{1}{3}xy$ could be preserved across all three pieces if the two cuts that form the sides of the middle piece are dependent on each other and coordinated such that the middle piece's total area is preserved, but this need not be true in general. In particular, this would require strictly that $b = c$ and $a = d$, implying $a + b = \frac{2}{3}y$ out of necessity that the middle portion retain base lengths of $\frac{1}{3}y$, which ought to always be true. This, of course, is a very particular construction of a diagonal cut, and it need not be true in general that $b = c$ and $a = d$. In fact, if this equality does not hold, equality is no longer preserved. Therefore, in general, this particular variant on the perturbed diagonal cutting method cannot yield perfect thirds of the original rectangular sandwich, except in the cases that $b = c$ and $a = d$, including the case that $a = b = c = d = \frac{1}{3}y$ . ◊

**Next Steps**

These results are an important step forward in sandwich mathematics. This further begs the question: What kinds of cuts *do* yield perfect thirds? What of curved cuts (curves that are regions of conic sections, exponentials, trigonometric functions, etc.)? Perhaps fractal cuts? Is there some family of cuts that share some invariant underlying structure? These methods of cutting, and perhaps consideration of optimizing cutting practices, warrant special attention that I hope to investigate in not-so-distant future studies.

**Conclusion and Consequences for Teaching and Learning Mathematics**

While this was admittedly a fun problem, a major asset in this approach is that the mathematics that are required are contained well within a sequence of secondary algebra courses. It is my contention that such a problem would be invaluable in the mathematics classroom for a number of reasons. The first is that this problem engages students in working deeply with systems of equations, constructing parametric relationships, and remaining aware of domain restrictions that are integral to considering the feasibility of their solutions. Another advantage is that this problem reveals the non-trivial nature of everyday activities. Too often in educational trajectories, we consider the reliance on 'real-world scenarios' to be comforting in that tangible activities provide satisfactory simplification. This exercise reveals that this need not be true, and that an activity that anybody can perform (in this case, cutting your favorite sandwich to share with two others) is not as straightforward as one might initially presume. Finally, this exercise reveals the promise that mathematics has in making sense of the order and/or disorder within our world. Teaching and learning mathematics ought to be a complex exchange of ideas and connections (O'Meara & Vaidya, 2021), and we in the mathematical community should encourage students to embrace all endeavors they might face. Working through this problem was incredibly fun and revealed a new truth about the world around us. For that, there will always be value in problem-solving.

In conclusion, diagonal cuts of a rectangular sandwich do not yield perfect thirds in general. While I am unsure what to do with this information, I am eager to bring this into all of my practices, regardless of which hat I might be wearing. ⋄